\documentclass[11pt]{amsart}
\usepackage{a4wide}
\usepackage{graphicx,enumitem}
\usepackage{amsmath,amssymb,amsthm}
\usepackage[english]{babel}
\usepackage{subcaption}
\usepackage{units}
\usepackage{diagbox}
\usepackage{fontenc}
\usepackage{xcolor}
\usepackage{hyperref}
\hypersetup{
    colorlinks=true,
    citecolor=green!50!black,
    linkcolor=red!80!black,
    filecolor=blue,
    urlcolor=blue!70!black,
}
\usepackage{caption}
\usepackage{csquotes}
\usepackage{ulem}
\newsavebox{\measurebox}
\usepackage{cancel}
\usepackage{tikz}
\usepackage{amssymb}
\DeclareSymbolFontAlphabet{\amsmathbb}{AMSb}%
\usepackage{mathbbol}
\usepackage{latexsym}
\usepackage{xspace,bm}
\usepackage{dsfont}
\usepackage[noabbrev]{cleveref}
\usepackage{cases}
\usepackage{todonotes}

\usepackage{mathtools}

\newcommand{\IP}{\amsmathbb{P}}
\newcommand{\R}{\amsmathbb{R}}

\newcommand{\C}{\amsmathbb{C}}

\newcommand{\IS}{\amsmathbb{S}}

\usepackage{float}

\newcommand{\cE}{\mathcal{E}}

\DeclareMathOperator{\E}{\amsmathbb{E}} %expectation
\newcommand{\dd}{\,\mathrm{d}}

\newtheorem{lemma}{Lemma}[section]
\newtheorem{proposition}[lemma]{Proposition}

\newtheorem{remark}[lemma]{Remark}

\usepackage{color}
\definecolor{darkgreen}{rgb}{0,.6,0}

\newcommand{\revision}[1]
{{\color{black} #1}}
\renewcommand{\Re}{\operatorname{Re}}
\renewcommand{\Im}{\operatorname{Im}}

\newcommand{\po}{{\psi_1}}
\newcommand{\pt}{{\psi_2}}
\title[Numerical signature of blow-up]{On the numerical signature of blow-up in matrix hydrodynamics}
\author[E.~Jansson]{Erik Jansson$^{*,\dagger}$}
\email[]{erikjans@chalmers.se}

\author[K.~Modin]{Klas Modin$^*$} 
\email[]{klas.modin@chalmers.se}

\address{$^*$Department of Mathematical Sciences, Chalmers University of Technology \& University of Gothenburg, S--412~96 G\"oteborg, Sweden.} 

\thanks{
$^\dagger$Corresponding author.
}

\begin{document}

\begin{abstract}
\revision{
Matrix hydrodynamics, namely spatial discretization via quantization theory, gives finite-dimensional models of hydrodynamic equations that keep their geometric structure.
We use it to study how the blow-up of a partial differential equation shows up in simulations. 
Our test case is the complexified Euler equations on the sphere, which are expected to be ill-posed in Sobolev spaces but well-posed for a short time for analytic data. 
We extend Zeitlin's model to these equations and linearize both the continuous and matrix equations around an imaginary rigid rotation. 
For the continuous linearization, random perturbations of Sobolev regularity leave the space of square-integrable functions instantly with probability one, whereas analytic perturbations do not up to an explicit time. 
The linearization of Zeitlin's model reproduces the continuous one up to truncation, with norm growth rates proportional to the matrix size. 
As Zeitlin's model, correctly discretized in time, quasi-preserves energy,  solutions in the matrix model remain bounded. The model itself does not blow up. 
The signature of blow-up instead how the time to reach a fixed relative size depends on the matrix size. 
For Sobolev data,  this time tends to zero as the matrix size grows, whereas for analytic data, its decrease slows down, consistent with the limit predicted by the linear theory. 
Simulations of the matrix model with an isospectral integrator are performed and agree with the linear predictions to within $0.1\%$ while conserving the spectrum to round-off. \\[1ex]
\textbf{Keywords:} Euler equations, blow-up, matrix hydrodynamics, Zeitlin's model, Euler-Arnold equations, PDE numerics, vorticity\\[1ex]
\textbf{MSC2010:} 35Q31, 35B44, 81S10, 65G50, 65M12}
\end{abstract}

\maketitle

\section{Introduction}\label{sec:intro}

Finite time blow-up in hydrodynamic equations is an active topic in the analysis of partial differential equations (PDEs).
For the 3-D incompressible Euler equations, significant results were given by Elgindi~\cite{El2021} and by Chen and Hou~\cite{ChHo2021}.
These results were partially guided by numerical simulations of Luo and Hou~\cite{Luo2014}, whose work thus exemplifies how careful numerical experiments can lead to theoretical breakthroughs.
On the other hand, in the description of the Millennium Prize Problem for the Navier--Stokes equations, Fefferman~\cite{FeffermanNS} writes: ``Many numerical computations appear to exhibit blow-up for solutions of the Euler equations, but the extreme numerical instability of the equations makes it very hard to draw reliable conclusions.''
Indeed, in numerical simulations of hydrodynamics, it is challenging to determine if the behavior you see persists in exact solutions or is a numerical artifact.

\revision{Numerical searches for singularities have a long history.
Early work include the study of the Taylor--Green vortex by Brachet et al. \cite{Brachet1983Smallscale} and the work on the 3D Euler equations by Kerr \cite{Kerr1993}. See for instance Gibbon \cite{Gibbon2008} or Hou \cite{Hou2009} for reviews. 

The difficulty is to tell finite-time blow-up from fast but finite growth.
For instance, for the surface quasi-geostropic equation, an apparent singularity suggested by Constantin et al.  \cite{Constantin1994} was later proven by C\'ordoba \cite{Cordoba1998} to be incompatible with the proposed blow-up scenario. 
Further, even for the 2D Euler equations, which are well-posed, the vorticity gradients can grow quickly, as shown by Denisov \cite{Denisov2009} and Zlato\v{s} \cite{Zlato2015}.
In the presence of a boundary, the growth can even be superexponential \cite{Kiselev2014}.
}

Matrix hydrodynamics  (cf.~\cite{MoVi2025b}) is the 
collective name for quantization-based discretizations of hydrodynamic-like equations. 
It began with the work of Zeitlin \cite{zeitlin1991, zeitlin2004}, for the 2-D Euler equations, based on quantization theory developed by Hoppe and others (\emph{cf.}~\cite{bordemann1991,bordemann1994,hoppe1989,le2018}).
Zeitlin's model captures all the geometric features of the 2-D Euler  equations. 

Whether blow-up occurs is tied to geometric properties in phase space.
For example, the global well-posedness of the 2-D Euler equations, contrasted against the (supposed) blow-up in the 3-D Euler equations, has to do with phase space geometry: in 2-D, vorticity is transported by the flow map as a function, which gives an abundance of conservation laws (Casimir function) not available for 3-D Euler.
Thus, due to the superior geometric properties, it is natural to ask how blow-up is manifested in matrix hydrodynamics.
\revision{Moreover, as a self-consistent but finite dimensional model, Zeitlin's method offers an additional challenge. Energy is conserved, and this bounds how large a solution can grow. Thus, to detect blow-up, a different signature than waiting for numerical overflow is needed.}

The approach is to carry out numerical experiments for an equation where we expect instant blow-up: the \emph{complexified Euler equations} \revision{for non-analytic initial data}, suggested by {\v{S}}verák~\cite{Sverak2018} as the natural complexification of Arnold's~\cite{arnold1966} geometric description of Euler's equations.
The equations are
\begin{equation}
    \label{eq:compeuler}
    \begin{gathered}
        \frac{\partial v}{\partial t} + \nabla_{\overline v} v+ \nabla_{v}^\top\, \overline{v}  = -\nabla p,\\
        \operatorname{div} v = 0,
    \end{gathered}
\end{equation}
where $v$ is a time-dependent complexified vector field, $p$ is a complex-valued function, $\overline{v}$ denotes the conjugated vector field, and $\nabla_{v}^\top$ denotes the adjoint of the complexified covariant derivative $\nabla_v$ relative to the standard $L^2$-pairing.

The interest in the complexified Euler equations \eqref{eq:compeuler} arises in their analysis, contrasted against the real setting.
Namely, the 2-D (real) Euler equations are globally well-posed, but {\v{S}}verák noticed that this global analysis does not extend to the complexified 2-D equations.
Indeed, 
Albritton and Ogden \cite{Albritton2023} showed that the complexified Euler equations~\eqref{eq:compeuler} on the doubly periodic square are ill-posed with instant blow-up (see details below). 
\revision{In this paper, we perform a linearized ill-posedness analysis in the spherical settings which indicates similar behaviour as in the doubly periodic square case.}
We thus have a candidate set of equations for studying how blow-up manifests itself in numerical simulations.

Just as small modifications of the equations can entirely change the analysis, one discretization to another can gravely affect the convergence behavior relative to blow-up.
In particular, for the 2-D complex Euler equations modified to the 2-D complex Navier-Stokes equations by adding a small amount of viscosity, it is expected that the standard well-posedness techniques for Navier--Stokes apply.
As {\v{S}}verák points out, ``This is no longer the case for the 2-D complex Euler equation, where the standard proofs of existence depend on more detailed properties of vorticity, which may not be shared by the complex equation.'' \cite{Sverak2018}.
The result of Albritton and Ogden verifies this statement and indeed makes it precise.
Consequently, to identify blow-up in the numerics, we better use a discretization that does not introduce artificial numerical viscosity, and matrix hydrodynamics is precisely such a discretization.

We now continue the introduction with a brief recollection of the standard Euler equations in 2-D and the basic idea behind Zeitlin's model.
Via the vorticity function $\omega=\operatorname{curl} v$, the Euler equations on the sphere can be written
\begin{equation}\label{eq:euler_2d_vorticity}
    \frac{\partial \omega}{\partial t} + \{\omega, \psi \} = 0, \quad -\Delta\psi = \omega,
\end{equation}
where $\{\cdot,\cdot\}$ is the Poisson bracket. 
The velocity field $v$ is recovered from the stream function $\psi$ by the skew gradient $v = \nabla^\top\psi$.
Based on the notion of quantization, where functions are replaced by operators (in this case matrices) and the Poisson bracket is replaced by the commutator, Zeitlin's model for the vorticity equation~\eqref{eq:euler_2d_vorticity} is the finite-dimensional isospectral matrix flow
\begin{equation}\label{eq:zeitlin}
    % \begin{gathered}
        \frac{dW}{dt} + \frac{1}{\hbar_N}[W, P] = 0, \quad
       \revision{ -\Delta_N P = W.}
    % \end{gathered}
\end{equation}
Here, the \emph{vorticity matrix} $W$ and the \emph{stream matrix} $P$ are  skew-Hermitian $N\times N$ matrices, and $\hbar_N = 2/\sqrt{N^2-1}$. 
The operator $\Delta_N$ is the \emph{Hoppe--Yau Laplacian} \cite{hoppe1998}, corresponding to the Casimir element for a finite-dimensional $\mathfrak{so}(3)$-representation, just as the Laplace--Beltrami operator on the sphere is the Casimir element for the infinite-dimensional representation of $\mathfrak{so}(3)$ on the space of smooth functions.
For more details on Zeitlin's model~\eqref{eq:zeitlin} and its geometric properties, see Zeitlin \cite{zeitlin1991,zeitlin2004}, Modin and Viviani \cite{modin2019,MoVi2022}, and Modin and Perrot \cite{MoPe2023}.
Gallagher \cite{Ga2002} established convergence of solutions of Zeitlin's model to solutions of the vorticity equation for the torus case originally considered by Zeitlin.
The convergence in the spherical case was established by Modin and Viviani~\cite{MoVi2024}.

\revision{In this paper, we linearise the complexified Euler equations on the sphere around the imaginary rigid rotation $\omega_0=iaz$.
The linearized system splits into subsystems for each degree $l$ and order $m$ of the spherical harmonics. 
We show that the linearized problem is ill-posed for data of finite Sobolev regularity, while for analytic data its solution exists up to a time $T^*$ set by the width of the strip of analyticity.
We show that Zeitlin's model reproduces the linearised dynamics exactly for all degrees up to $N-1$, and that conservation of energy bounds its solutions by a constant that grows with the matrix size $N$.
From this we derive bounds on the time at which the expected norm of a random perturbation reaches a fixed threshold.
For data of finite Sobolev regularity this time tends to zero, at {\color{black}least as fast as } $\log N/N$, while for analytic data it tends to a positive limit below $T^*$.
In numerical experiments up to $N=512$, the threshold times of the full model agree with those of the linear solution to within $0.1\%$.
The dependence of the threshold time on $N$ is therefore a signature that separates instant blow-up from a solution that exists for a positive time.
}

\revision{
The paper is organized as follows.
In \cref{sec:geomprelim}, we provide some geometric preliminaries needed for the later presentation.
In \cref{sec:complexeuler}, we present the complexified Euler equations and the linear analysis around the imaginary rigid rotation.
In \cref{sec:zeitlin}, we derive Zeitlin's model for the complexified equations, show that it reproduces the linearised dynamics exactly, and derive the predicted threshold times.
In \cref{sec:experiments}, we compare these predictions with numerical experiments and validate the computations.
We end with a discussion in \cref{sec:discussion}.
}

\subsection*{Acknowledgments}
This work was supported by the Swedish Research Council (grant number 2022-03453), the Knut and Alice Wallenberg Foundation (grant number KAW~2024.0229), the Göran Gustafsson Foundation for Research in Natural Sciences and Medicine, and the Wallenberg AI, Autonomous Systems and Software Program (WASP) funded by the Knut and Alice Wallenberg Foundation.
The computations were enabled by resources provided by the National Academic Infrastructure for Supercomputing in Sweden (NAISS), partially funded by the Swedish Research Council through grant agreement no. 2022-06725.

\section{Geometric preliminaries}
\label{sec:geomprelim}
In this section, we describe the geometric setting underlying the 2-D Euler equations on the sphere.
First, the sphere $\IS^2\subset\R^3$ is a K\"ahler manifold.
\revision{For $x\in\IS^2$ and $u,v\in T_x\IS^2$, its metric $g$ is the Euclidean inner product $g_x(u,v)=u\cdot v$, its symplectic form is the area form $\Omega_x(u,v)=x\cdot(u\times v)$, and its complex structure is $J_xv=x\times v$.
Thus $J$ rotates each tangent vector by $\pi/2$, counterclockwise as seen from outside the sphere.
These satisfy $J^2=-\mathrm{id}$ and $g(\cdot,\cdot)=-\Omega(J\cdot,\cdot)$.}

The gradient $\nabla f$ of a smooth function $f \in C^\infty(\IS^2,\R)$ is defined by 
\begin{align*}
    \dd f(\cdot) =  g(\nabla f,\cdot),
\end{align*}
where $\dd$ is the exterior derivative. 
The Hamiltonian vector field $X_f$ of $f$ is similarly defined by the relation
\begin{align*}
    \dd f(\cdot) =\Omega(X_f,\cdot).  
\end{align*}
Since 
\begin{align*}
    \dd f(\cdot) = g(\nabla f,\cdot) = -\Omega(J \nabla f,\cdot),
\end{align*}
we have that $-J\nabla f = X_f $. 
This means that Hamiltonian vector fields are obtained by point-wise $\pi/2$ rotation of corresponding gradient vector fields.

Quantization unfolds from the Poisson algebra of smooth functions $(C^\infty(\IS^2,\R),\{\cdot,\cdot\})$, where the Poisson bracket $\{\cdot,\cdot\}\colon C^\infty(\IS^2,\R)\times C^\infty(\IS^2,\R) \to C^\infty(\IS^2,\R)$ is defined by
\begin{align}
    \label{eq:poi}
    \{f,g\} = \Omega(X_f,X_g) = X_g(f). 
\end{align}
In the standard, unit radius embedding $\IS^2 \subset \R^3$, it is explicitly given by
\[
    \{f,g\}(\mathbf{r}) %= (\nabla f(\mathbf{r})\times \mathbf{r})\cdot \nabla g(\mathbf r) 
    = (\mathbf{r}\times \nabla g(\mathbf{r}))\cdot \nabla f(\mathbf r).
\]
We sometimes also use spherical coordinates, with the convention $$(x_1,x_2,x_3) = (\sin\theta\cos\phi,\sin\theta \sin\phi,\cos\theta),$$ where $(\theta,\phi) \in (0,\pi)\times (0,2\pi)$. In these coordinates, away from the poles, $\nabla f = (\partial_\theta f, \partial_\phi f/\sin\theta)$.

The Poisson bracket is skew-symmetric and satisfies Leibniz's rule and the Jacobi identity. 
In particular, it is a Lie bracket. 
Further,
\begin{align}
    \label{eq:poi-iso}
    [X_f,X_g]_{\mathfrak{X}} = -X_{\{f,g\}},
\end{align}
where $[\cdot,\cdot]_{\mathfrak{X}}$ is the Lie bracket of vector fields. 
\Cref{eq:poi-iso} means that $f\mapsto X_f$ is a Lie algebra anti-morphism between $(C^\infty(\IS^2,\R),\{\cdot,\cdot\})$ and $(\mathfrak{X}_{\Omega}(\IS^2),[\cdot,\cdot]_{\mathfrak{X}})$, where $\mathfrak{X}_{\Omega}(\IS^2)$ denotes Hamiltonian vector fields on $\IS^2$. 
The space of all vector fields is denoted $\mathfrak{X}(\IS^2)$.
Thus, $\mathfrak{X}_{\Omega}(\IS^2)$ is a subalgebra of $\mathfrak{X}(\IS^2)$.

\section{The complexified Euler equations}
\label{sec:complexeuler}
The complexified Euler equations arise naturally as the \emph{Euler--Arnold equations} (cf.~\cite{ArKh2021}) on the complexification of the Lie algebra of Hamiltonian vector fields, as derived by {\v{S}}verák~\cite{Sverak2018}.

In general, if $\mathfrak{g}$ is a real Lie algebra, finite or infinite dimensional, its complexification $\mathfrak{g}_\C = \mathfrak{g} \otimes \C$ consists of elements of the form $X + iY$, where $X,Y \in \mathfrak{g}$.
The corresponding complexified Lie bracket is thus given by 
\begin{equation}
	[X_1+iY_1,X_2+iY_2]_\C = [X_1,X_2] - [Y_1,Y_2]+ i[X_1,Y_2] + i[X_2,Y_1].
\end{equation}
\revision{Furthermore, an inner product $(\cdot\,,\cdot)$ on $\mathfrak{g}$ extends to a Hermitian inner product on $\mathfrak{g}_\C$, linear in the first and conjugate-linear in the second argument.
Its real part is an inner product on $\mathfrak{g}_\C$ regarded as a real vector space.}
The complexification of $\mathfrak{X}_\Omega(\IS^2)$ is $\mathfrak{X}_\Omega(\IS^2)_\C = \mathfrak{X}
_\Omega(\IS^2) \otimes \C$.  
The $L^2$-inner product on $\mathfrak{X}_\Omega(\IS^2)$ extends to the Hermitian inner product, 
\begin{align}
	(v_1,v_2)_{L^2} = \int_{\IS^2} v_1 \cdot \overline{v_2},
\end{align}
where $v_1, v_2 \in \mathfrak{X}_\Omega(\IS^2)_\C$.
\revision{The complexified Euler equations are the Euler--Arnold equations \eqref{eq:compeuler} for the real part of this inner product, with $\mathfrak{X}_\Omega(\IS^2)_\C$ regarded as a real Lie algebra.}
As such, {\v{S}}verák's~\cite{Sverak2018} derivation is given in terms of the complex vector field $v$.
Next, we give instead a derivation based on complexification of the vorticity formulation~\eqref{eq:euler_2d_vorticity}.

\subsection{ The complexified vorticity equation}

Just as there is a vorticity formulation \eqref{eq:euler_2d_vorticity} of the standard 2-D Euler equations, we can derive a vorticity formulation of the 2-D complexified Euler equations~\eqref{eq:compeuler}. 
It may be obtained by applying $\operatorname{curl}$ to \cref{eq:compeuler}, but more geometrically, we derive it as an Euler--Arnold equation on the complexified symplectic vector fields $\mathfrak{X}_\Omega(\IS^2)\otimes\C$. 
First, the complex vector field $v$ is given as the sum of two real Hamiltonian vector fields, which we see as the Hamiltonian vector field of a complex stream function.
Indeed, the Hamiltonian vector field of a purely imaginary function $if \in iC^\infty(\IS^2,\R)$ is given by 
\begin{align*}
    X_{if} = iX_f. 
\end{align*}
Thus, given a complex-valued stream function $\psi \in C^\infty(\IS^2,\C)$, the corresponding complex Hamiltonian vector field is
\begin{align}
    \label{eq:decompos}
    X_\psi = X_{\Re \psi + i \Im \psi } =X_{\Re \psi} + iX_{\Im \psi}.
\end{align}

Next, the Poisson algebra structure on $C^{\infty}(\IS^2,\C)$ is readily obtained by complexifying $(C^\infty(\IS^2,\R),\{\cdot,\cdot\})$.
Explicitly, the bracket on $C^{\infty}(\IS^2,\C)$ is given by  
\begin{align*}
        &\{\po,\pt\}_\C = \{\Re \po,\Re \pt\}-\{\Im \po,\Im \pt\}\\
        &\qquad+ i \left(\{\Im \po,\Re \pt\}+\{\Re \po,\Im \pt\}\right).
\end{align*}
 By properties inherited from the real-valued case, 
     $(C^\infty(\IS^2,\C),\{\cdot,\cdot\}_\C)$ is a Poisson algebra. 
    
The infinitesimal action of a Hamiltonian vector field on an imaginary function is determined by the Poisson bracket.
Indeed, we have that
\begin{align*}
     X_{ig}(f) = i\{f,g\},
 \end{align*}
     where $f,g \in C^\infty(\IS^2,\R)$. 
 Thus, for $\po,\pt \in C^\infty(\IS^2,\C)$, 
 \begin{align}
     \label{eq:poidef}
     \begin{split}
         &X_\pt(\po) = X_{\Re \pt}(\po) +iX_{\Im \pt} (\po) \\
         &= X_{\Re \pt}(\Re \po)+i X_{\Re \pt}(\Im \po) + X_{i\Im \pt}(\Re \po)+iX_{i\Im \pt}(\Im \po)  \\
%         &=\{\Re \po,\Re \pt\}-\{\Im \po,\Im \pt\}\\
        &\qquad +i\left(\{\Re \po,\Im\pt\}+\{\Im \po,\Re \pt \}\right) =\{\po,\pt\}_{\C}.
    \end{split}
\end{align}

By the definition of the real Lie bracket of vector fields and \cref{eq:poidef}, one shows that the Lie bracket of complexified vector fields is related to the complex Poisson bracket by
\begin{align}
    \label{eq:pamorph}
    -X_{\{\po,\pt\}_\C} = [X_\po,X_\pt]_{\mathfrak{X}\otimes \C}.  
\end{align}

\Cref{eq:pamorph} does not define an isomorphism. 
Indeed, all constant functions are mapped to the zero vector field.
Since constant functions form the center of the Poisson algebra $(C^\infty(\IS^2,\C),\{\cdot,\cdot\}_\C)$, we obtain that the reduced Poisson algebra $(C^\infty(\IS^2,\C)/\C,\{\cdot,\cdot\}_\C)$ \revision{is a Lie algebra anti-isomorphic} to $\mathfrak{X}_\Omega(\IS^2)_\C$.

On a general Lie algebra $\mathfrak{g}$, the Euler--Arnold equations are
\begin{align}\label{eq:abstract}
     \frac{d\mu}{dt} = \operatorname{ad}_{dH(\mu)}^* \mu,
\end{align}
where $\mu \in \mathfrak{g}^*$ and $H \in C^\infty(\mathfrak{g}^*,\R)$ is the Hamiltonian function with exterior derivative $d H(\mu) \in \mathfrak{g}^{**}\simeq \mathfrak{g}$.  
The mapping $\operatorname{ad}^*\colon\mathfrak{g} \to \operatorname{End}(\mathfrak{g}^*)$ is the dual, under the canonical dual pairing, of the adjoint representation $\operatorname{ad}\colon\mathfrak{g} \to \operatorname{End}(\mathfrak{g})$. 
This way, \cref{eq:compeuler} is retrieved with $\mathfrak{g} = \mathfrak{X}_\Omega(\IS^2)_\C$ and the Hamiltonian 
\begin{align}
    \label{eq:ham_vec}
    H(v) = \frac{1}{2}(v,v)_{L^2}.
\end{align}

As mentioned, the mapping $\psi \mapsto X_\psi$ provides an anti-isomorphism between the two Lie algebras $(C^\infty(\IS^2,\C)/\C,\{\cdot,\cdot\}_\C)$ and $\mathfrak{X}_\Omega(\IS^2)_\C$, which yields the following.

\begin{proposition}
    The vorticity formulation of 
\cref{eq:compeuler} is, 
\begin{align}
    \label{eq:vort}
    % \begin{gathered}
        \frac{\partial\omega}{\partial t} = \{-\overline{\psi},\omega\}_{\C},\quad
         -\Delta \psi =  \omega.
    % \end{gathered}
\end{align} 
\end{proposition}

\begin{proof}
    The first step towards a complexified vorticity formulation of \cref{eq:compeuler} is to derive the corresponding Euler--Arnold equation on $(C^\infty(\IS^2,\C)/\C,\{\cdot,\cdot\}_\C)$.
     As the Euler--Arnold equation evolves on the dual of the algebra, we note that the (smooth) dual of $(C^\infty(\IS^2,\C)/\C,\{\cdot,\cdot\}_\C)$
is the annihilator of constant functions in the space of top-forms. 
A top-form $\beta = \omega \Omega$, where $\omega \in C^\infty(\IS^2,\C)$, acts on a function $\psi \in C^\infty(\IS^2,\C)$ by 
\begin{align}
    \label{eq:l2pair}
    \revision{\langle \beta,\psi\rangle = \Re\int_{\IS^2} \omega \bar{\psi} .}
\end{align}
Since top-forms can be identified with smooth functions, the (smooth) dual of $C^\infty(\IS^2,\C)/\C$ is isomorphic to $C^\infty_0(\IS^2,\C)$, the smooth complex-valued functions with vanishing mean. 

The next step is to compute $\operatorname{ad}_\psi^*\colon C^\infty_0(\IS^2,\C) \to C^\infty_0(\IS^2,\C)$. 
By definition,
\begin{align*}
    &\langle  \operatorname{ad}_\psi^* \omega, \xi\rangle =  
    \langle   \omega, \operatorname{ad}_\psi \xi\rangle = 
    \Re\int_{\IS^2} \omega \overline{\{\psi,\xi\}_{\C}}   = 
    \Re\int_{\IS^2} \omega\{\bar\psi,\bar\xi\}_{\C}   = \\
    &\Re\int_{\IS^2} \left( \{\bar\psi,\omega\bar\xi\}_{\C} - \bar\xi\{\bar\psi,\omega\}_{\C} \right)  = 
     \langle \{-\bar\psi,\omega\}_\C, \xi\rangle,
\end{align*}
where, in the fourth equality, Leibniz's rule is used. 
Therefore, 
\begin{equation}
    \label{eq:adstar_def}
    \operatorname{ad}_\psi^* \omega = \{-\bar\psi,\omega\}_\C\revision{.}
\end{equation}

The final step is to find the canonical relationship between $\omega$ and $\psi$ determined by the Hamiltonian of the complexified Euler equation, i.e.,  \cref{eq:ham_vec}. 

The Hamiltonian $\tilde{H}: C^\infty_0(\IS^2,\C) \to \R$ is given by 
\begin{align}
    \label{eq:stepHam}
    \tilde H(\omega) =  \frac{1}{2}\int_{\IS^2}  X_\psi \cdot \bar X_\psi.
\end{align}
Since 
\begin{align*}
     \int_{\IS^2}  X_\psi \cdot \bar X_\psi = \int_{\IS^2}  \nabla  \psi \cdot \nabla \overline{\psi} ,
\end{align*} 
integration by parts yields that \cref{eq:stepHam} becomes 
\begin{align}
    \begin{split}
         \tilde H(\omega) = \frac{1}{2}\int_{\IS^2} -\Delta  \psi  \overline{\psi},
    \end{split}
    \label{eq:Ham}
\end{align}
so $\tilde H(\omega) = \frac{1}{2}(-\Delta  \psi,\psi )_{L^2}$, and we have that 
\begin{equation}
        \label{eq:hamdef}
      d \tilde H(\omega) = -\Delta^{-1} \omega,
\end{equation}
Therefore, by inserting \cref{eq:adstar_def,eq:hamdef} into the abstract formulation 
\eqref{eq:abstract}, we see that the vorticity equation is given by \cref{eq:vort}.
\end{proof}

\subsection{Properties of the complexified Euler equations}

We first remark that there are several conserved quantities. 
First, note that the Hamiltonian~\eqref{eq:ham_vec} (i.e., the kinetic energy)
is conserved. 
Furthermore, the complex vorticity formulation \eqref{eq:vort} reveals an infinite set of Casimir functions: if $f\colon \C \to \C$ is a holomorphic function, then 
\begin{equation}\label{eq:casimirs}
    C_f(\omega) = \int_{\IS^2} f(\omega)
\end{equation}
is conserved.
Indeed,
\[
    \frac{d}{dt} C_{f}(\omega) = \int_{\IS^2} f'(\omega)\{-\overline{\psi},\omega \}_{\C}
    = \int_{\IS^2} \overline{\psi} \{ f'(\omega),\omega\}_{\C}  ,
\]
which vanishes since $f'$ can be expanded in a Taylor series and $\{\omega^k,\omega\}_{\C}=0$.
The Casimirs \eqref{eq:casimirs} reveal a special structure of the complexified Euler equations in two dimensions, analogous to the special structure of the real Euler equations in two dimensions.

Moreover, we note that if the initial value is real, then the vorticity remains real and the complexified vorticity equations~\eqref{eq:vort} reduce to the standard vorticity equations~\eqref{eq:euler_2d_vorticity}.
Indeed, \cref{eq:vort} separates neatly into components: 
with $\omega = \omega_R + i \omega_I$ and $\psi = \psi_R + i \psi_I$ it can be written
\begin{align*}
    \frac{\partial \omega_R}{\partial t}  = -\{\psi_R,\omega_R\} -\{\psi_I,\omega_I\} \\
    \frac{\partial \omega_I}{\partial t} = -\{\psi_R,\omega_I\}+\{\psi_I,\omega_R\}.
\end{align*}

\revision{
\subsection{Linear ill-posedness of an imaginary rigid rotation}
In order to know which flows to study numerically and which regularity of perturbations to use, we employ, as in Albritton and Ogden \cite{Albritton2023} and {\v{S}}verák~\cite{Sverak2018}, linear ill-posedness. 
We consider perturbations of the imaginary rigid rotation 
\begin{align*}
    \omega_0 = iaz, 
\end{align*}
where $a \in \R^+$ and $z = \cos(\theta)$ with $\theta \in [0,\pi)$ is the polar angle on the sphere. 
Note that $\omega_0$ is proportional to the spherical harmonic $Y_{1,0}$, so it is zonal and is stationary, as $\Delta^{-1} Y_{1,0} = -\frac{1}{2}Y_{1,0}$ so that the Poisson bracket in \cref{eq:vort} vanishes at $\omega_0$. 
Consider a perturbation $f_0$ of $\omega_0$. 
Using that $\omega_0$ is zonal, the linearised equation for $f(t)$, evolving from $f_0$ is given by 
\begin{align}
\label{eq:lin}
\frac{\partial f}{\partial t} = ia \left(\partial_\phi f/2 + \partial_\phi \overline{\psi_f}\right), ~f(0) = f_0. 
\end{align}
where $\psi_f = \Delta^{-1}f$.
We show that this system is ill-posed. 
To do this, we add random perturbations, either $H^s(\IS^2)$-valued random fields (see e.g. \cite{LS15}), where $s \in [0,\infty)$, 
or analytic noise, for which the spectrum decays exponentially in the $\ell$ wavenumber.

\begin{proposition}
    \label{prop:lin_growth}
    Let $f_0$ be a random perturbation of $\omega_0$ of the form 
    \begin{align*}
    \varepsilon\sum_{l=1}^\infty \sum_{m=-l}^l c_{l,m}Y_{l,m}
    \end{align*}
    where $\varepsilon>0$ and the real and imaginary parts of $c_{l,m}$ are independent normal random variables with mean zero and variance $\sigma_l^2$.
    If the perturbation is slowly decaying, that is, if  $\lim_{l\to\infty}\frac{-\log\sigma_l}{l}=0,$ 
    then  the solution $f(t)$ to \cref{eq:lin} satifies 
    \begin{align*}
        \|f(t)\|_{L^2(\IS^2)} = \infty,
    \end{align*}
    almost surely for all $t>0$.
    If instead the noise is analytic, that is, 
     if $(\sigma_l^2\lesssim e^{-2\delta l})$  for some $\delta>0$, then there is a $T^*$ such that the expected $L^2(\IS^2)$-norm of $f(t)$ is finite for $0\le t< T^*$ and infinite for $t\geq T^*$, where $T^*=2\delta/|a|$   
\end{proposition}
\begin{proof}
    We expand $f(t)$ in spherical harmonics.
We note that $\overline{Y_{l,m}}=(-1)^mY_{l,-m}$, so we group the $m$ and $-m$ terms to obtain, 
\begin{align*}
f=\sum_{l\ge1}f_{l,0}Y_{l,0}+\sum_{m\ge1}\sum_{l\ge m}\left(u_{l,m}Y_{l,m}+\overline{v_{l,m}Y_{l,m}}\right),
\end{align*}
where $v_{l,m}=(-1)^m\overline{u_{l,-m}}$.
Inserting into the linearised equation and comparing the coefficients of $Y_{l,m}$ and of $\overline{Y_{l,m}}$ gives the system $\dot U = A_{l,m}U$ where $U = (u_{l,m},v_{l,m})^T$ and
\begin{align}
    \label{eq:coeff}
    A_{l,m} = am\begin{pmatrix}
-1/2 & -\mu_l\\
\mu_l & 1/2
\end{pmatrix},
\end{align}
where $\mu_l = 1/(l(l+1))$. Since $A_{l,m}^2 = a^2 m^2 \left(1/4-\mu_l^2\right)I$, the solution to \cref{eq:coeff}
is 
\begin{align*}
     U(t) = e^{tA_{l,m}}U(0) = \left(\cosh\left(t\rho_{l,m}\right)I + V_{l,m}A_{l,m}\right) U(0),
\end{align*}
where $\rho_{l,m} = |a|m\sqrt{1/4-\mu_l^2}$ and $V_{l,m} = \sinh(t\rho_{l,m})/\rho_{l,m}$ if $\rho_{l,m} \neq 0$ and $V_{l,m} = t$ if $\rho_{l,m} = 0$.
The matrix $A_{1,1}$ is not zero, but $A_{1,1}^2=0$, so $e^{tA_{1,1}}=I+tA_{1,1}$ grows only linearly in $t$.
If $l\ge2$, then $\mu_l \le1/6$, so $\rho_{l,m}\ge|a|m\sqrt{1/4-1/36}=\tfrac{\sqrt2}{3}|a|m>0$, and the eigenvalues of $A_{l,m}$ are $\pm\rho_{l,m}$.
Moreover $\rho_{l,m}\le|a|m/2$, and $\rho_{l,m}$ increases with both $l$ and $m$.

The trace of $A_{l,m}$ is zero and the sum of the squares of its entries is $a^2m^2\left(\frac12+2\mu_l^2\right)$.
Hence the sum of the squares of the entries of $e^{tA_{l,m}}$ is $2B_{l,m}(t)$, where
\begin{align*}
B_{l,m}(t)=\cosh^2(t\rho_{l,m})+a^2m^2\left(\frac14+\mu_l^2\right)V_{l,m}^2.
\end{align*}
For $l\ge2$, the second term equals $\kappa_l\sinh^2(t\rho_{l,m})$ with $\kappa_l=\frac{\frac14+\mu_l^2}{\frac14-\mu_l^2}\le\frac54$, so $B_{l,m}(t)\le\frac94\cosh^2(t\rho_{l,m})\le\frac94e^{|a|mt}$.

We have that $u_{l,m}(0)=\varepsilon c_{l,m}$ and $v_{l,m}(0)=(-1)^m\varepsilon\,\overline{c_{l,-m}}$.
Both are complex normal variables whose real and imaginary parts are independent with variance $\varepsilon^2\sigma_l^2$, so $\E|u_{l,m}(0)|^2=\E|v_{l,m}(0)|^2=2\varepsilon^2\sigma_l^2$.

By independence of the coefficients, for every $t\ge0$ we have that
\begin{align}
\label{eq:elin}
\E\|f(t)\|_{L^2(\IS^2)}^2=2\varepsilon^2\sum_{l\ge1}\sigma_l^2+4\varepsilon^2\sum_{m\ge1}\sum_{l\ge m}\sigma_l^2B_{l,m}(t).
\end{align}

Now let $m\ge2$ and consider the pair $l=m$, for which $\rho_{m,m}>0$.
 $A_{m,m}$ is real with the real eigenvalue $\rho_{m,m}$. 
 The corresponding eigenvector $w_m$ satisfies
  $w_m^Te^{tA_{m,m}}=e^{t\rho_{m,m}}w_m^T$.
Write $U_{m,m}(t)=(u_{m,m}(t),v_{m,m}(t))^T$, so that $U_{m,m}(t)=e^{tA_{m,m}}U_{m,m}(0)$.
By the Cauchy--Schwarz inequality,
\begin{align*}
|U_{m,m}(t)|\ge|w_m^TU_{m,m}(t)|=e^{t\rho_{m,m}}|w_m^TU_{m,m}(0)|.
\end{align*}
Since $w_m$ is a real unit vector, $w_m^TU_{m,m}(0)$ is a complex normal variable whose real and imaginary parts are independent with variance $\varepsilon^2\sigma_m^2$, so 
\begin{align*}
\IP\left(E_m \right)=e^{-1/2}.
\end{align*}
Where $E_m$ is the event that $|w_m^TU_{m,m}(0)|^2\ge\varepsilon^2\sigma_m^2$.

On $E_m$ we have, for every $t\ge0$,
\begin{align*}
|u_{m,m}(t)|^2+|v_{m,m}(t)|^2\ge\varepsilon^2\sigma_m^2e^{2t\rho_{m,m}}.
\end{align*}
Due to the independence of the coefficients, we note that the sum of probabilities is divergent, so by the second Borel--Cantelli lemma, with probability one infinitely many of them occur.

For slowly decaying noise this gives the following:
Since $m\ge2$, $\rho_{m,m}\ge\frac{\sqrt2}{3}|a|m$, so on $E_m$
\begin{align*}
|u_{m,m}(t)|^2+|v_{m,m}(t)|^2\ge\varepsilon^2\left(\sigma_me^{\frac{\sqrt2}{3}|a|tm}\right)^2,
\end{align*}
and the right-hand side tends to infinity as $m\to\infty$, for every $t>0$.
Hence, with probability one, $\|f(t)\|_{L^2(\IS^2)}=\infty$ for all $t>0$.
In words, noise of any Sobolev regularity leaves $L^2$ at once under the linearised flow.

For analytic noise with rate $\delta$, let $T^*=2\delta/|a|$.
For $0\le t<T^*$ the expected norm is finite.
Indeed, the term with $l=1$ is finite, and for $l\ge2$ we have $\sigma_l^2B_{l,m}(t)\le\frac94e^{-2\delta l}e^{|a|mt}$, so the sum over $l\ge m$ is at most $\frac94e^{(|a|t-2\delta)m}/(1-e^{-2\delta})$, and the sum over $m$ converges when $|a|t<2\delta$.
So $f(t)\in L^2$ with probability one for each $t<T^*$.
For $t\ge T^*$, on $E_m$ we have $|u_{m,m}(t)|^2+|v_{m,m}(t)|^2\ge\varepsilon^2e^{2t\rho_{m,m}-2\delta m}$.
Since $2t\rho_{m,m}\ge2T^*\rho_{m,m}=2\delta m\sqrt{1-4/(m(m+1))^2}$ and $\sqrt{1-x}\ge1-x$ for $0\le x\le1$, the exponent is at least $-8\delta/m^3\ge-\delta$.
Therefore, with probability one, infinitely many terms are at least $\varepsilon^2e^{-\delta}$, and $\|f(t)\|_{L^2(\IS^2)}=\infty$ for all $t\ge T^*$.

\end{proof}

}

\section{Zeitlin's model for the complexified Euler equations}

\label{sec:quant}\label{sec:zeitlin}
In this section, we show how a finite-dimensional analogue of the developments in the previous section naturally leads to a generalization of Zeitlin's approach to the vorticity formulation~\eqref{eq:vort} of the complexified Euler equation on $\IS^2$.

We first review how the Poisson algebra $(C^\infty(\IS^2,\C),\{\cdot,\cdot \})$ is discretized via quantization.
Further details are found in the papers by Hoppe and Yau~\cite{hoppe1998} and by Modin and Viviani~\cite{modin2019,MoVi2024}.

The vorticity equation evolves on the (smooth) dual of the space of symplectic vector fields. 
The quantized equation should have an equivalent structure: it should evolve on the dual of a matrix Lie algebra.
To obtain the quantization scheme, we use a matrix analog of the Laplacian, as given by Hoppe and Yau~\cite{hoppe1998}.
It is an operator $\Delta_N\colon \mathfrak{sl}(N,\C)\to \mathfrak{sl}(N,\C)$ which allows us to build an eigenbasis of $\mathfrak{sl}(N,\C)$ corresponding to the spherical harmonic basis of functions on $\IS^2$.

Explicitly, we begin with a spin $\frac{N-1}{2}$ representation of  $\mathfrak{so}(3)$ in $\mathfrak{u}(N)$. 
Its generators $X_1,X_2,X_3$ satisfy the commutation relations of $\mathfrak{so}(3)$ up to scaling by $\hbar_N = 2/\sqrt{N^2-1}$:
\begin{align*}
[X_2,X_1] = \hbar_N X_3,\quad [X_1,X_3] = \hbar_N X_2,\quad [X_3,X_2] = \hbar_N X_1.
\end{align*}
The infinite-dimensional analogue is the following. 
Via the natural embedding $\IS^2 \subset \R^3$, there is a representation of $\mathfrak{so}(3)$ in $C^\infty(\IS^2,\R)$, with generators given by the coordinate functions $x_1,x_2,x_3$. 
They each generate rotations about their respective coordinate axis. 
This representation satisfies the  commutation relations of $\mathfrak{so}(3)$,
\begin{align*}
\{x_2,x_1\} = x_3,\quad \{x_1,x_3\} = x_2,\quad \{x_3,x_2\} = x_1,
\end{align*}

Moreover, the Laplace--Beltrami operator on $\IS^2$ can be expressed in terms of the Poisson bracket and the generators as  
\begin{align}
    \label{eq:lap}
    \Delta = \sum_{k=1}^3 \{x_k,\{x_k,\cdot\}\}.  
\end{align}
Analogous to this expression for $\Delta$, the Hoppe--Yau Laplacian on $\mathfrak{sl}(N,\C)$ is defined by
\begin{align}
    \label{eq:HY}
    \Delta_N = \frac{1}{\hbar_N^2} \sum_{k=1}^3 [X_k,[X_k,\cdot]].
\end{align}

The eigenvalues of the Hoppe--Yau Laplacian $\Delta_N$ coincide with the first $N^2$ eigenvalues of the Laplace--Beltrami operator $\Delta$, see \cite{hoppe1998}. 
The quantization of a function $\zeta\in C^\infty(\IS^2,\C)$ is now given by a mapping from $C^\infty(\IS^2,\C)$ to $\mathfrak{gl}(N,\C)$ as follows.

The expansion of  $\zeta$ in  the spherical harmonics $\{Y_{l,m}\mid l\geq 0, \, m = -l,\ldots,l\}$  is truncated, 
    \begin{align*}
    \zeta^N  = \sum_{l=0}^N \sum_{m = -l}^l c_{l,m} Y_{l,m}.
\end{align*} 
The quantization is then performed by identifying the first $N^2$ spherical harmonics, being eigenfunctions of $\Delta$, with corresponding eigen-matrices of $\Delta_N$, denoted $T_{l,m}^N$.
Once such a correspondence is established, we obtain a mapping from functions to matrices via
\begin{align*}
    \zeta^N \mapsto Z = \sum_{l=0}^N \sum_{m = -l}^l c_{l,m}  T_{l,m}^N. 
\end{align*}

It remains to describe the matrices $T_{l,m}^N$. 
The quantized Laplacian $\Delta_N$ maps $\pm m$-diagonal matrices to $\pm m$-diagonal matrices. 
A unique eigenbasis to $\Delta_N$ is given if we take $T_{l,m}$ to be $\pm m$-diagonal. 
Under this restriction, the eigenvalue problem 
\begin{align*}
    \Delta_N T_{l,m}^N = -l(l+1) \, T_{l,m}^N 
\end{align*}
has a unique solution up to sign {\color{black} and scale.} 
We determine the sign based on the corresponding continuous harmonics $Y_{l,m}$.
{\color{black}We fix the scale so that 
\begin{equation*}
    \| T_{l,m}^N\|_{2}^2 \coloneqq \frac{4\pi}{N}\operatorname{tr}\big((T_{l,m}^N)^* T_{l,m}^N\big) = 1.
\end{equation*}
Thus, this scaled Frobenius norm corresponds to the $L^2$ norm on $C^\infty(\IS^2,\C)$.
Furthermore, as part of the framework of quantization, the spectral norm $\|W \|_\infty \coloneqq \sup_{|\boldsymbol{v}|=1}|W\boldsymbol{v}|$ corresponds to the $L^\infty$ norm (cf.~\cite[Table~1]{MoVi2024}).
}
This process provides a quantization of complex-valued functions.
From the point of view of complexification, just as $C^\infty(\IS^2,\C) \simeq C^\infty(\IS^2,\R)\otimes \C$, the quantized equivalent of complex-valued functions is $\mathfrak{gl}(N,\C) \simeq \mathfrak{u}(N)\otimes \C$. 

In direct analogue to the previous section, we also obtain the $\operatorname{ad}_P^*$ operator on $\mathfrak{sl}(N,\C)$ via
\[
    \langle\operatorname{ad}^*_{P}W, X\rangle = \Re\operatorname{tr}(W[P,X]^*) = \Re\operatorname{tr}(-W[P^*,X^*])=
   \Re \operatorname{tr}([P^*,W]X^*) = \langle [P^*,W], X\rangle .
\]
Thus, for $P\in \mathfrak{sl}(N,\C)$ and $W\in \mathfrak{sl}(N,\C)$, $\operatorname{ad}^*_PW = [P^*,W]$, which is well-defined on $\mathfrak{sl}(N,\C)$, since $\C^* I$ is the center of $\mathfrak{gl}(N,\C)$.
Notice that classical conjugation $\psi\mapsto \overline\psi$ corresponds in the quantized framework to $P\mapsto -P^*$.

\subsection{Matrix vorticity equations}
We are now in position to state the spatially discretized version of the vorticity formulation \eqref{eq:vort} of the complexified Euler equations on $\IS^2$.
Indeed, the equations are
\begin{align}
\label{eq:quantvort}
    % \begin{gathered}
        \frac{d W}{dt} = \frac{1}{\hbar_N}[P^*,W],\qquad
        -\Delta_N P =  W,
    % \end{gathered}
\end{align}
where the isomorphism
$-\Delta_N \colon \mathfrak{sl}(N,\C) \to \mathfrak{sl}(N,\C)$ is the Hoppe-Yau Laplacian above.
{\color{black}The corresponding Hamiltonian is
\begin{equation*}
    \tilde H_N(W) = \frac{2\pi}{N}\operatorname{tr}(P^* W).
\end{equation*}
}

Notice that an initially skew-Hermitian solution remains skew-Hermitian.
Indeed, \cref{eq:quantvort} separate into skew-Hermitian and Hermitian components, just as the complexified vorticity equation~\eqref{eq:vort} separates into real and imaginary components.

We remark that the matrix vorticity equations~\eqref{eq:quantvort} are still not fully discrete: a time discretization is needed, i.e., a numerical integrator. 
To ensure that the geometric structure of Euler-Arnold equations are properly reflected in the numerics, it is advisable to use an integrator that preserves the Lie--Poisson structure.
See, for instance, Eng{\o} and Faltinsen \cite{engo2001}, Bogfjellmo and Marthinsen \cite{Bogfjellmo2015}, and Modin and Viviani \cite{modin2019_1} for examples of such integrators.

\begin{remark}
    The matrix vorticity equations \eqref{eq:quantvort} preserve the structure that leads to the Casimirs~\eqref{eq:casimirs}.
    Indeed, if $f$ is holomorphic, then
    \begin{equation}\label{eq:zeitlin_casimirs}
        C^N_f(W) = {\color{black}\frac{4\pi}{N}}\operatorname{tr}(f(W))
    \end{equation}
    is a Casimir.
    These Casimirs reflect that the flow is isospectral.
\end{remark}
\revision{
Note that while the complexified Euler equations \eqref{eq:compeuler} may exhibit blow-up, this is not possible in a finite-dimensional and energy preserving method. 
Indeed, we note that 
% \begin{align*}
% \|W\|_2^2\le\|W\|_F^2=C^2\sum_{l,m}|W_{l,m}|^2\le N(N-1) C^2\sum_{l,m}\frac{|W_{l,m}|^2}{l(l+1)}=2N(N-1)E,
% \end{align*}
{\color{black}
\begin{multline*}
    \|W\|_\infty^2\le \operatorname{tr}(W^*W)=\frac{N}{4\pi}\|W \|_2 = \frac{N}{4\pi}\sum_{l,m}|W_{l,m}|^2 \le \\ \le \frac{N^2(N-1)}{4\pi} \sum_{l,m}\frac{|W_{l,m}|^2}{l(l+1)}=\frac{N^2(N-1)}{2\pi}\tilde H_N(W),
\end{multline*}
}
% \todo[inline]{I think it is better to introduce the energy functional $E = E(W)$ right after the complexified Euler equations are introduced. Now it's implicitly defined in the proof.}
% where $\|\cdot\|_F$ denotes the Frobenius norm, $\|\cdot\|_2$ is the spectral norm and $C$ is the normalization constant of the matrix spherical harmonics.
This means that the energy {\color{black}determined by the initial conditions $\tilde H_N(W_0)$} bounds on how large the matrix can get. 
Typically, the solution grows until it reaches the ceiling, where is starts to oscillate. 
To capture the numerically capture blow up, we should be able to replicate theoretical predictions of how $N$ determines the norm growth, that is, how quickly the norm reaches some prescribed value. 
To do this, we replicate the linearization argument from \cref{sec:quant}.
The Zeitlin equivalent of the data $\omega_0 = ia z$ is given by 
\begin{equation}
W_0=\tfrac{ia}{C}T_{1,0},
\end{equation}
where $C$ is a normalization constant giving unit norm to the matrix spherical harmonics. 
By Modin and Viviani \cite{modin2019}, the matrix $T_{1,0}$ is diagonal and $W_0$ satifies 
\begin{equation}
\tfrac1\hbar_N[W_0,T_{l,m}]=-am\,T_{l,m},\qquad\tfrac1\hbar_N[W_0,T_{l,m}^*]=am\,T_{l,m}^*.
\end{equation} 

For an initial perturbation $F_0$, the linearized Zeitlin system is given by 
$\dot F=\mathcal{M}F$, with $\mathcal{M}F=\tfrac1\hbar_N\left([P_0^*,F]+[P_F^*,W_0]\right)$ and $-\Delta_NP_F=F$.
By repeating the arguments in \cref{sec:complexeuler}, we have that the equation for the coefficients of $F(t)$ is 
\begin{equation}
\begin{pmatrix}\dot u_{l,m}\\\dot v_{l,m}\end{pmatrix}=am\begin{pmatrix}-1/2&-\mu_l\\\mu_l&1/2\end{pmatrix}\begin{pmatrix}u_{l,m}\\v_{l,m}\end{pmatrix},
\end{equation}
which is the same matrix $A_{l,m}$ as for the continuous equation, so Zeitlin's model reproduces the linear dynamics of the continuous equation at this background exactly until the truncation $N$.
As a consequence, the growth rates of the linearized Zeitlin equations are the same as that of the linearized continuous equations, 
and the largest growth rate is
\begin{equation}
R_N=\rho_{N-1,N-1}=|a|(N-1)\sqrt{\tfrac14-\tfrac1{(N-1)^2N^2}}.
\end{equation}

We can now make some concrete predictions for growth of the matrices.
\begin{proposition}
    Consider a perturbation of the form $F(0)=\varepsilon\sum_{l=1}^{N-1}\sum_{m= -l}^l c_{l,m}T_{l,m}$, with $c_{l,m}$ as in \cref{prop:lin_growth} and  $\sigma_l$ independent of $N$.
    For $\theta>0$, let the predicted time $T_N(\theta)$ be the first time at which $\E[\|F(t)\|_F]\ge\theta$.
    Then, 
    \begin{enumerate}[nosep,label=(\roman*)]
        \item $T_N(\theta)$ is non-increasing in $N$.
        \item For slowly decaying noise, $T_N(\theta)\to0$ as $N\to\infty$.
        \item For analytic noise with rate $\delta$, $T_N(\theta)\to T_\infty(\theta)$ as $N\to\infty$, where $T_\infty(\theta)$ is the unique time with $\mathbb E\left[\|f(t)\|^2\right]=\theta^2$, and $0<T_\infty(\theta)<T^*=2\delta/|a|$, where $f(t)$ is the solution to \cref{eq:lin}.
    \end{enumerate}
    Moreover, if  $\sigma_l \propto l^{-p}$, then $T_N(\theta)$ is at most of order $\log N/N$.
For $p>1$,  $T_N(\theta)$ is at least of order $1/N$.
\end{proposition}
\begin{proof}
%The formulas below are stated for $\tau=1$.
%For general $\tau$, every norm of $F$ is multiplied by $\tau$, so a threshold $\Theta$ for $\|F\|_F$ corresponds to the threshold $\theta=\Theta/\tau$ below.

%In the runs, $T_{\rm lin}$ is compared with the time $T(K)$ at which $\|W(t)-W_0\|_F\ge K\|W_0\|_F$.
%There, $\tau=\sqrt N$ and $\|W_0\|_F=\frac{|a|}{c}\,\tau=10\sqrt N$, so $\theta=K\|W_0\|_F/\tau=10K$.

For every $t\ge0$, as the Zeitlin coefficient system is exactly the same as the undiscretized one,
\begin{align}
    \label{eq:nN}
\E\left[\|F(t)\|_F^2\right]=2\varepsilon^2\sum_{l=1}^{N-1}\sigma_l^2+4\varepsilon^2\sum_{m=1}^{N-1}\sum_{l=m}^{N-1}\sigma_l^2\,B_{l,m}(t),
\end{align}
but with the sums stopped at $l=N-1$.
We note that the norm is non-decreasing in $t$ and tends to the infinite-dimensional value as $N\to\infty$.
Recall the infinite-dimensional value $\E\|f(t)\|_{L^2(\IS^2)}^2$.
It may be infinite. 
Let $\theta > n_\infty(0)$ and let let $T_N(\theta)$ be the first time at which 
\begin{align*}
\E\left[\|F(t)\|_F^2\right]\geq\theta^2.
\end{align*}

Write $n_N(t)^2=\E\|F_{\rm lin}(t)\|_F^2/C^2$ and $n_\infty(t)^2=\E\|f(t)\|^2$, which may be infinite.
Each $B_{l,m}$ is continuous and non-decreasing in $t$, and increasing $N$ adds non-negative terms, so $n_N$ is continuous and non-decreasing in $t$, non-decreasing in $N$, and $n_N(t)\to n_\infty(t)$ by monotone convergence.
This gives (i).

For (ii), fix $t>0$.
For $m\ge2$, the terms of \eqref{eq:elin} with $l=m$ satisfy $4\varepsilon^2\sigma_m^2B_{m,m}(t)\ge\varepsilon^2\sigma_m^2e^{2t\rho_{m,m}}\ge\varepsilon^2\big(\sigma_me^{\frac{\sqrt2}3|a|tm}\big)^2$, since $B_{m,m}(t)\ge\cosh^2(t\rho_{m,m})\ge\frac14e^{2t\rho_{m,m}}$.
By the assumption on $\sigma_m$, these terms tend to infinity, so $n_\infty(t)=\infty$.
Since $n_N(t)$ increases to $n_\infty(t)$, $n_N(t)\ge\theta$ for all large $N$, that is, $T_N(\theta)\le t$.
Since $t>0$ was arbitrary, $T_N(\theta)\to0$.

For (iii), $n_\infty$ is finite on $[0,T^*)$ by \cref{prop:lin_growth}.
It is strictly increasing on $[0,T^*)$, because the terms with $l\ge2$ are.
At $t=T^*$, the terms with $l=m$ are at least $\varepsilon^2e^{-\delta}$ by the estimate in the analytic case, so $n_\infty(T^*)=\infty$, and since $n_\infty(t)$ is strictly increasing, $n_\infty(t)\to\infty$ as $t\uparrow T^*$.
Hence there is a unique $T_\infty(\theta)\in(0,T^*)$ with $n_\infty(T_\infty(\theta))=\theta$.
Since $n_N\le n_\infty$, $T_N(\theta)\ge T_\infty(\theta)$.
For $T_\infty(\theta)<t<T^*$, $n_\infty(t)>\theta$, so $n_N(t)\ge\theta$ for all large $N$, that is, $T_N(\theta)\le t$.
Hence $T_N(\theta)\to T_\infty(\theta)$.

Let $\bar R_N=\frac{|a|(N-1)}2$, so that $R_N=\bar R_N\sqrt{1-4/((N-1)^2N^2)}\le\bar R_N$.
For $l\ge2$, $B_{l,m}(t)\le\frac{9}{4}e^{|a|mt}\le\frac94e^{2\bar R_Nt}$, and $B_{1,1}(t)=1+\frac{a^2t^2}2\le e^{|a|t}\le\frac94e^{2\bar R_Nt}$.
Every term of \eqref{eq:nN} is therefore at most $\frac94e^{2\bar R_Nt}$ times its value at $t=0$, so $n_N(t)\le\frac32\,n_N(0)\,e^{\bar R_Nt}$.
Keeping only the term with $l=m=N-1$, and using $B_{N-1,N-1}(t)\ge\cosh^2(R_Nt)\ge\frac14e^{2R_Nt}$, gives $n_N(t)\ge\varepsilon\,\sigma_{N-1}\,e^{R_Nt}$.
Hence
\begin{equation}\label{eq:TNbounds}
\frac1{\bar R_N}\log\frac{2\theta}{3\,n_N(0)}\le T_N(\theta)\le\frac1{R_N}\log\frac{\theta}{\varepsilon\,\sigma_{N-1}}.
\end{equation}
For $\sigma_l=C_\sigma l^{-p}$,  the upper bound in \eqref{eq:TNbounds} gives $T_N(\theta)\le\big(p\log(N-1)+\log\frac{\theta}{\varepsilon C_\sigma}\big)/R_N$, so $T_N(\theta)$ is at most of order $\log N/N$.
For $p>1$, $n_N(0)\le n_\infty(0)<\infty$, and the lower bound shows that $T_N(\theta)$ is at least of order $1/N$.

\end{proof}
}

\section{Numerical study of blow-up}
\label{sec:experiments}
We now numerically study the blow-up of the complexified Euler equations on $\IS^2$. 
To this end, we first remark that any purely real solution (in the Zeitlin case, purely skew-Hermitian) to the \revision{complexified Euler equations devolves to a solution of the real Euler equations}. 
Therefore, we focus on perturbations of purely imaginary stationary solutions (in the matrix case, purely Hermitian).
Indeed, the blow-up of the complexified Euler equations arises from instability of the interaction of the real and imaginary components.

In the numerical study, we use the matrix vorticity equations~\eqref{eq:quantvort}.
\revision{Following \cref{sec:zeitlin}, we consider perturbations of $W_0 = (ia/C) T_{1,0}$, with $ a/C = 10$. 
The skew-Hermitian part is therefore zero, and the initial flow is zonal. }
Zonal flows are a natural choice for the initial condition, as zonal flows are stationary solutions to the Euler equations.
Further, to construct the perturbations, consider the random matrix given by 
\begin{align*}
   \revision{ \widetilde{\cE}_0 = \varepsilon \sum_{l=1}^{N-1} \sum_{m = -l}^l c_{l,m}  T_{l,m}^N,}
\end{align*}
\revision{Here,  $\varepsilon$ is a small real number, in all simulations chosen to be $0.1$, and $c_{l,m}$ are random complex numbers chosen such that the real and imaginary parts are independent normal random variables with mean $0$ and variance $\sigma_l^2$, i.e., we use one variance per eigenspace. 
The final perturbation $\cE_0$ is then given by normalizing $\widetilde{\cE}_0$ to unit Frobenius norm.}

In all Zeitlin runs, we use the isospectral midpoint method of \cite{modin2019} to \cref{eq:quantvort} to numerically integrate the equations in time.
The method is a Lie--Poisson integrator that preserves the Lie--Poisson structure of the equations, and in particular, the Casimirs. 
Each step of the integration is given by      \revision{
\begin{align}
    \label{eq:isomp}
    \begin{split} 
        &W_n = \left(I-\frac{\delta t}{2\hbar_N} (\Delta_N^{-1} \widetilde W)^*\right) \widetilde{W} \left(I+\frac{\delta t}{2\hbar_N} (\Delta_N^{-1} \widetilde W)^*\right), \\
        &W_{n+1} = \left(I+\frac{\delta t}{2\hbar_N} (\Delta_N^{-1} \widetilde W)^*\right) \widetilde{W} \left(I-\frac{\delta t}{2\hbar_N} (\Delta_N^{-1} \widetilde W)^*\right),
    \end{split}
\end{align}}
where $I$ is the identity matrix and $\widetilde{W}$ is an intermediate variable obtained by solving the first implicit step. 
\revision{In the simulations, we compute $\widetilde{W}$ at each step by fixed-point iterations, converged to a tolerance of $10^{-13}$. Unless otherwise specified, we use a time step size of $\delta t = 10^{-5}\,\hbar_N$. 
Throughout, let $T(K)$ denote the first time that $\|W(t)-W_0\|_{F} \geq K \|W_0\|$.

We run four experiments and one set of validation runs. 
\subsection{Growth and saturation }
\label{sec:expA}
We first illustrate how the solutions grow until a fixed cap, depending on $N$ and as predicted in \cref{sec:zeitlin}. 
We use $\sigma_l = (1+l(l+1))^{-3/2}$ to generate the perturbations $\cE_0$. 
We are this in the setting of slowly decaying noise. 
We run the simulation for $N = 32, 64, 128, 192, 384$ and $512$ until the final time $T_{\mathrm{end}} = 6/N$. 
In \cref{fig:expA} we plot the relative deviation from the  Frobenius norm of $\|W_0\|_F$,that is 
\begin{align*}
    D(W(t))\coloneqq \|W(t)-W_0\|_F/\|W_0\|_F.
\end{align*}

\begin{figure}[h!]
    \centering
    \includegraphics[width=0.75\linewidth]{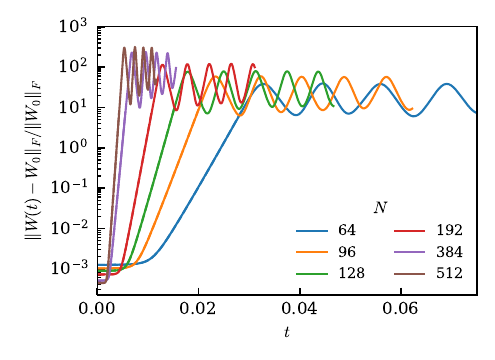}
    \caption{Growth of the relative deviation from the initial value in the Frobenius norm. Note how for all values $N$, the norm grows until it saturates at some level. However, the growth rate increases as $N$ increases.}
    \label{fig:expA}
\end{figure}

We observe the expected behavior, in line with the predictions in \cref{eq:compeuler}. 
The relative deviation $D(W(t))$ rises until some ceiling,  after which it oscillates.
As $N$ increases, the rate of growth of $D(W(t))$ increases.  
We show the rate in \Cref{fig:expA_rate}, and observe a good agreement with a linear growth of the rate, in line with the predictions in \cref{sec:zeitlin}.
This indicates that the growth rate is unbounded as $N \to \infty$, which is the discrete counterpart of the ill-posedness of the linearised equation.

\begin{figure}[ht!]
    \centering
    \includegraphics[width=0.75\linewidth]{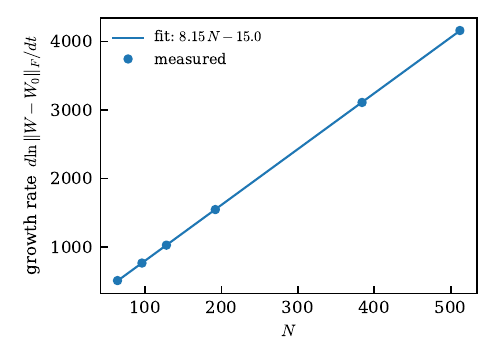}
    \caption{Growth rate of $D(W(t))$ at $t=T(1)$ against $N$, with a linear fit.}
    \label{fig:expA_rate}
\end{figure}

\subsection{Growth rate behavior: slowly decaying case}
\label{sec:expB}
In this experiment, we study how the rate of growth behaves as a function of $N$, when the noise is slowly decaying. 
In this setting, the prediction is that the rate of growth should increase with $N$. 
In other words, $T(K)$, the time it takes for the relative deviation to become larger than $1$,  should decrease.
The linear theory in \cref{sec:zeitlin} predicts that for values of $N$ considered here, the time should fall off as at most $\log(N)N^{-1}$, indicating that as $N \to \infty$, the growth is immediate. 

We consider four different noise classes of the form $(1+l(l+1))^{-\alpha}$, with $\alpha = 0, 1/2, 1$ and $-3/2$. 
Each run is otherwise structured as in \cref{sec:expA}, and we re-generate the perturbation and re-run for a total of $15$ times at each $\alpha$ and $N$. 
The average results are plotted against $N$ in \cref{fig:expB}. 
We also plot the predicted linear predictions, by computing the linearized solution $\cref{eq:lin}$ starting from the same initial value at the same truncation value. 
\begin{figure}[h!]
    \centering
    \includegraphics[width=0.75\linewidth]{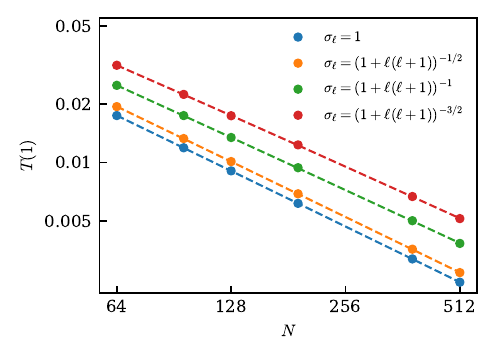}
    \caption{The time it takes to reach a fixed relative deviation $D(W(t))$ against the matrix size $N$. Zeitlin predictions are dots, and the linear prediction is the dashed line. Note the close-to linear decay and the striking agreement with the linear preduction.}
    \label{fig:expB}
\end{figure}
We note a striking agreement between Zeitlin's method and the linearized prediction. In fact, the agreement is to with in $0.1\%$.
While the linear ill-posedness in principle does not say anything about the nonlinear system, the agreement indicates that linear prediction might have some validity. 

\subsection{Growth rate behavior: analytic case}
\label{sec:expC}
The linear theory predicted that in the analytic case, we should see that that the time to reach some fixed cut-off approaches an $N$-independent final time upper-bounded by $T^* =2\delta/|a|$, where $\delta$ is the radius of analyticity, or, the decay rate in $\sigma_l = e^{-\delta l}$. 
We repeat the same set-up as in \cref{sec:expB} for $\delta = 0.025$ and $\delta = 0.05$.  In  \cref{fig:expC} we plot the average  observed times along with the corresponding linear predictions extended to high matrix sizes. 
Along the analytic cases, the highest regularity run from \cref{sec:expB} is included. 
\begin{figure}[h!]
    \centering
    \includegraphics[width=0.75\linewidth]{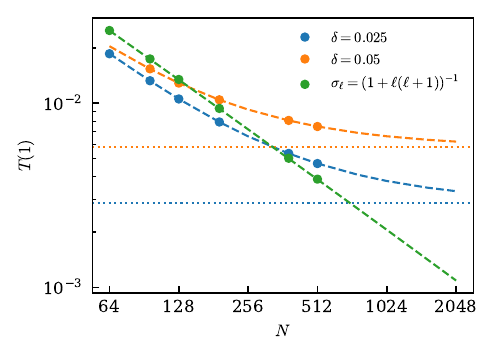}
    \caption{The time it takes to reach a fixed relative deviation $D(W(t))$ against the matrix size $N$ for different perturbations. Zeitlin predictions are dots, and the linear prediction is the dashed line. The linear predictions are extended to higher values of $N$. Note the close-to linear decay and the striking agreement with the linear prediction.}
    \label{fig:expC}
\end{figure}
We observe in \cref{fig:expC} that the analytic perturbations lead to a qualitatively different behavior. 
The observed times level off and move towards the predicted threshold, as indicated by the linear predictions. 
We have omitted the high-$N$ Zeitlin runs, as the computational improvements required for such sizes at date require skew-Hermitian settings and as such do not cover the present setting.

\subsection{Validation of the experimental set-up}
First, to validate that the method behaves correctly, we validate the structure-preserving qualities of the method, namely, that the eigenvalue drift is neglible, that the relative energy error is independent of $N$, and that the system is reversible, i.e., that running until saturation and then reversing the flow recovers the initial condition. 
We first re-use the runs from \cref{sec:expA} and compute the maximum relative energy error $|E(W(t)-E(W(0))|/E(W(0))$, the maximum over the runs of the maximum eigenvalue drift $\max_{j}|\lambda_j(W(t)-\lambda_j(W(0))|$, where $\lambda_j(A)$ is the $j$th eigenvalue of the matrix $A$. 
To compute the recovery error, for each run, the setup of \cref{sec:expA} us run backwards using negative step size for the same length of time $T_{\mathrm{end}}$, producing the solution $W^{-}(T_{\mathrm{end}})$. 
As the system \eqref{eq:zeitlin} is reversible, and the used method is reversible, the recovery error
\begin{align*}
   \| W^{-}(T_{\mathrm{end}})-W(0)\|/\|W(0)\|
\end{align*}
should be small. 
We report the results in \cref{tab:methodworks}. 
\begin{table}
\label{tab:methodworks}
\centering
\begin{tabular}{r c c c}
\hline
$N$ & max.\ rel.\ energy error & max.\ eigenvalue drift & recovery error \\
\hline
64 & $3.1\times10^{-9}$ & $3.0\times10^{-11}$ & $2.4\times10^{-8}$ \\
96 & $3.8\times10^{-9}$ & $3.3\times10^{-11}$ & $3.4\times10^{-6}$ \\
128 & $3.4\times10^{-9}$ & $6.7\times10^{-11}$ & $5.1\times10^{-6}$ \\
192 & $3.5\times10^{-9}$ & $9.2\times10^{-11}$ & $3.3\times10^{-6}$ \\
384 & $3.8\times10^{-9}$ & $2.7\times10^{-10}$ & $1.2\times10^{-5}$ \\
512 & $3.9\times10^{-9}$ & $7.8\times10^{-10}$ & $2.2\times10^{-5}$ \\
\hline
\end{tabular}
\caption{Validation of the structure-preserving properties of the method on the runs from \cref{fig:expA}.}
\end{table}
We observe that all quantities behave as expected. 
It is also important to validate that the chosen step size is reasonable, in that the the times $T(K)$ can be resolved. 
To do this, we fix $N = 128$ and let the stepsize vary between  $\delta t = 10^{-3}\,\hbar_N, \delta t = 10^{-4}\,\hbar_N, \ldots, \delta t = 10^{-6}\,\hbar_N$. 
We use $\sigma_l = (1+(l+1)l)^{-3/2}$ and perform $5$ repeats. In \cref{tab:stepsize}, we report the median time to reach $T(K)$ for $K = 1,2,4$ and $8$.   
\begin{table}
    \centering
    \begin{tabular}{l c c c c}
\hline
$\delta t/\hbar_N$ & $T(1)$ & $T(2)$ & $T(4)$ & $T(8)$ \\
\hline
$10^{-6}$ & 0.01344 & 0.01411 & 0.01478 & 0.01545 \\
$10^{-5}$ & 0.01345 & 0.01413 & 0.01478 & 0.01545 \\
$10^{-4}$ & 0.01359 & 0.01422 & 0.01484 & 0.01547 \\
$10^{-3}$ & 0.01406 & 0.01563 & 0.01563 & 0.01563 \\
\hline
\end{tabular}
    \caption{Validation of the chosen stepsize. We observe that the times to reach various relative deviations settle at $\delta_t = 10^{-5}\hbar_N$. }
    \label{tab:stepsize}
\end{table}
We observe that the times to reach various relative deviations settle at $\delta_t = 10^{-5}\hbar_N$.

}

\revision{
\section{Discussion}\label{sec:discussion}
A finite-dimensional model that preserves (or in case of the temporally discretized system, quasi-preserves) energy cannot blow up, as the initial energy restricts how large the norms can become. 
However, as the dimension increases, the maximum attainable norm increases with $N$. Therefore, to detect blow-up with a structure preserving model we cannot hope to wait for the norm of the solution to exceed some predetermined large cut-off. Instead, we use the time $T(K)$ at which the deviation from the background reaches $K$ times the norm of the background, and study how it depends on $N$.

For slowly decaying random perturbations of a zonal initial condition, the linear theory predicts predicts that this time tends to zero, at most of order $\log N/N$. 
For analytic data, it predicts a positive limit below $T^*=2\delta/|a|$.

The simulations follow both predictions, and the threshold times of the full model agree with those of the linear solution to within $0.1\%$.
The decrease of $T(K)$ with $N$ is therefore a signature of instant blow-up, and its levelling off indicates a solution that exists for a positive time.

We remark that the analysis is linear. 
The agreement with the nonlinear simulations indicate that the threshold is reached while the dynamics still behave linearly.
In this case, the signature does not depend on the nonlinear phase.
}
\bibliographystyle{acm} %%%% .BST file

\bibliography{references} %%%%% .Bib file

\end{document}